\documentclass[11pt, letterpaper, leqno]{amsart}
\usepackage{amsmath}
\usepackage{amsthm}
\usepackage{amsrefs}
\usepackage{lipsum}
\usepackage{amsopn}
\usepackage{amssymb}
\usepackage{fullpage}
\numberwithin{equation}{section}
\usepackage{float}
\usepackage{graphicx}
\usepackage{titlesec}
\usepackage{fancyhdr}
\usepackage{fancyhdr}
\newtheorem{question}{Question}[section]

\newtheorem{theorem}{Theorem}[section]

\theoremstyle{definition}

\titleformat{\section}
  {\normalfont\scshape}{\thesection}{1em}{\centering}
\titleformat{\subsection}
{\normalfont\scshape}{\thesubsection}{1em}{\centering}

\begin{document}
\title{Hyperbolic distance and membership of conformal maps in the Hardy space}

\author{Christina Karafyllia}  
\address{Department of Mathematics, Aristotle University of Thessaloniki, 54124, Thessaloniki, Greece}
\email{karafyllc@math.auth.gr}   
\thanks{I thank Professor D. Betsakos, my thesis advisor, for his advice during the preparation of this work, and the Onassis Foundation for the scholarship I receive during my Ph.D. studies.}

\fancyhf{}
\renewcommand{\headrulewidth}{0pt}

\fancyhead[RO,LE]{\small \thepage}
\fancyhead[CE]{\small Hyperbolic distance and membership of conformal maps in Hardy space}
\fancyhead[CO]{\small Christina Karafyllia} 
\fancyfoot[L,R,C]{}



\begin{abstract} Let $\psi$ be a conformal map of the unit disk $\mathbb{D}$ onto an unbounded domain and, for $\alpha >0$, let ${F_\alpha }=\left\{ {z \in \mathbb{D}:\left| {\psi \left( z \right)} \right| = \alpha } \right\}$. If ${H^p}\left( \mathbb{D} \right)$ denotes the classical Hardy space and $d_\mathbb{D} {\left( {0,{F_\alpha }} \right)}$ denotes the hyperbolic distance between $0$ and $F_\alpha$ in $\mathbb{D}$, we prove that $\psi$ belongs to ${H^p}\left( \mathbb{D} \right)$ if and only if
	\[\int_0^{ + \infty } {{\alpha ^{p - 1}}{e^{ - {d_{\mathbb{D}}}\left( {0,{F_\alpha }} \right)}}d\alpha }  <  + \infty .\]
This result answers a question posed by P. Poggi-Corradini.
\end{abstract}

\maketitle
\section{Introduction}\label{section1}

A classical problem in geometric function theory is to find the Hardy number of a region by looking at its geometric properties (see e.g. \cite{Ba}, \cite{Ha} and \cite{Ra}). Answering a question of P. Poggi-Corradini (\cite[p. 36]{Co}), we give a necessary and sufficient integral condition for whether a conformal map of $\mathbb{D}$ belongs to ${H^p}\left( \mathbb{D} \right)$ by studying the hyperbolic metric in its image region. For a domain $D$, a point $z \in D$ and a Borel subset $E$ of $\overline D $, let ${\omega _D}\left( {z,E} \right)$ denote the harmonic measure at $z$ of $\overline E$ with respect to the component of $D \backslash \overline{E}$ containing $z$. The function ${\omega _D}\left( { \cdot ,E} \right)$ is exactly the solution of the generalized Dirichlet problem with boundary data $\varphi  = {1_E}$   (see \cite[ch. 3]{Ahl} and \cite[ch. 1]{Gar}). The hyperbolic distance between two points $z,w$ in the unit disk $\mathbb{D}$ (see \cite[ch. 1]{Ahl}, \cite[p. 11-28]{Bea}) is defined by 
\[{d_\mathbb{D}}\left( {z,w} \right) = \log \frac{{1 + \left| {\frac{{z - w}}{{1 - z\bar w}}} \right|}}{{1 - \left| {\frac{{z - w}}{{1 - z\bar w}}} \right|}}\]
and the Green function for $\mathbb{D}$ (see \cite[p. 41]{Gar}) is directly related to the hyperbolic distance, as it is obvious by its definition
\[{g_\mathbb{D}}\left( {z,w } \right) = \log \left| {\frac{{1 - z\bar w }}{{z - w }}} \right|= \log \frac{{{e^{{d_\mathbb{D}}\left( {z,w } \right)}} + 1}}{{{e^{{d_\mathbb{D}}\left( {z,w} \right)}} - 1}}.\]
They are both conformally invariant and thus they can be also defined on any simply connected domain $D \ne \mathbb{C}$ by considering a conformal map $f$ of $\mathbb{D}$ onto $D$. Then ${d_D}\left( {z,w} \right) = {d_\mathbb{D}}\left( {{f^{ - 1}}\left( z \right),{f^{ - 1}}\left( w \right)} \right)$ and ${g_D}\left( {z,w} \right) = {g_\mathbb{D}}\left( {{f^{ - 1}}\left( z \right),{f^{ - 1}}\left( w \right)} \right)$ for every $z,w \in D$. Therefore, for $z,w \in D$,

\begin{equation}\label{gree}
{g_D}\left( {z,w } \right) = \log \frac{{{e^{{d_D}\left( {z,w} \right)}} + 1}}{{{e^{{d_D}\left( {z,w} \right)}} - 1}}.
\end{equation}
Also, for a set $E \subset D$, we define ${d_D}\left( {z,E} \right): = \inf \left\{ {{d_D}\left( {z,w} \right):w \in E} \right\}$.

The Hardy space with exponent $p$, $0<p<+\infty$, and norm ${\left\|  \cdot  \right\|_p}$ (see \cite[p. 1-2]{Du}, \cite[p. 435-441]{Gar}) is defined to be
\[{H^p}\left( \mathbb{D} \right) = \left\{ {f \in H\left( \mathbb{D} \right):\left\| f \right\|_p^p = \mathop {\sup }\limits_{0 < r < 1} \int_0^{2\pi } {{{\left| {f\left( {r{e^{i\theta }}} \right)} \right|}^p}d\theta  <  + \infty } } \right\},\]
where $H\left( \mathbb{D} \right)$ denotes the family of all holomorphic functions on $\mathbb{D}$. S. Yamashita \cite{Ya} proved that a function $f \in H\left( \mathbb{D} \right)$ belongs to ${H^p}\left( \mathbb{D} \right)$, $0<p<+\infty$, if and only if 
\[ \int_\mathbb{D} {{{\left| {f\left( z \right)} \right|}^{p - 2}}{{\left| {f'\left( z \right)} \right|}^2}\log \frac{1}{{\left| z \right|}}dA\left( z \right)}  <  + \infty, \]
where $dA$ is the Lebesgue measure on $\mathbb{D}$. The fact that a function $f$ belongs to ${H^p}\left( \mathbb{D} \right)$ imposes a restriction on the growth of $f$ and this restriction is stronger as $p$ increases. If $f$ is a conformal map on $\mathbb{D}$, then $f \in {H^p}\left( \mathbb{D} \right)$ for all $p<1/2$ (\cite[p. 50]{Du}). 

Harmonic measure and hyperbolic distance are both conformally invariant and several Euclidean estimates are known for them. Thus, expressing the ${H^p}\left( \mathbb{D} \right)$-norms of a conformal map $\psi$ on $\mathbb{D}$ in terms of harmonic measure and hyperbolic distance, we are able to obtain information about the growth of the function by looking at the geometry of its image region $\psi \left( {\mathbb{D}} \right)$. Indeed, if $\psi $ is a conformal map on $\mathbb{D}$ and ${F_\alpha } = \left\{ {z \in \mathbb{D}:\left| {\psi \left( z \right)} \right| = \alpha } \right\}$ for $\alpha >0$, then P. Poggi-Corradini proved  (see \cite[p. 33]{Co}) that
\begin{equation}\label{isod}
\psi  \in {H^p}\left( \mathbb{D} \right) \Leftrightarrow \int_0^{ + \infty } {{\alpha ^{p - 1}}{\omega _{\mathbb{D}}}\left( {0,{F_\alpha }} \right)d\alpha }  <  + \infty.
\end{equation}

He also proved that the Beurling-Nevanlinna projection theorem (see \cite[p. 43-44]{Ahl}, \cite[p. 9-10, 35]{Co}) implies that

\begin{equation}\label{1.1}
{e^{ - {d_\mathbb{D}}\left( {0,{F_\alpha }} \right)}}\le \frac{\pi}{2 }{\omega _\mathbb{D}}\left( {0,{F_\alpha }} \right) .
\end{equation}
This observation led him to state two questions (see \cite[p. 36]{Co}):

\begin{question}\label{quest} Let $\psi $ be a conformal map  of $\mathbb{D}$ onto an unbounded domain and,  for $\alpha >0$, let ${F_\alpha } = \left\{ {z \in \mathbb{D}:\left| {\psi \left( z \right)} \right| = \alpha } \right\}$. Does there exist a positive constant $K$ such that for every $\alpha >0$,
	\[{\omega _\mathbb{D}}\left( {0,{F_\alpha }} \right) \le K{e^{ - {d_\mathbb{D}}\left( {0,{F_\alpha }} \right)}}?\]
\end{question}

\begin{question}\label{que} More generally, is it true that
	\[\psi  \in {H^p}\left( {\mathbb{D}} \right) \Leftrightarrow \int_0^{ + \infty } {{\alpha ^{p - 1}}{e^{ - {d_{\mathbb{D}}}\left( {0,{F_\alpha }} \right)}}d\alpha }  <  + \infty ?\]	
\end{question}

We gave a negative answer to the first one in \cite{Ka} and the following theorem provides a positive answer to the second question.

\begin{theorem}\label{theo} Let $\psi $ be a conformal map of $\mathbb{D}$ onto an unbounded simply connected domain $D$ and ${F_\alpha } = \left\{ {z \in \mathbb{D}:\left| {\psi \left( z \right)} \right| = \alpha } \right\}$  for $\alpha >0$. If $0<p<+\infty$ then
	\[\psi  \in {H^p}\left( {\mathbb{D}} \right) \Leftrightarrow \int_0^{ + \infty } {{\alpha ^{p - 1}}{e^{ - {d_{\mathbb{D}}}\left( {0,{F_\alpha }} \right)}}d\alpha }  <  + \infty .\]	
\end{theorem}

\section{Proof of Theorem \ref{theo}}\label{section}

\proof Suppose that $\psi  \in {H^p}\left( {\mathbb{D}} \right)$ for some $0<p<+\infty$. This in conjunction with (\ref{isod}) and (\ref{1.1}) implies that
\[\int_0^{ + \infty } {{\alpha ^{p - 1}}{e^{ - {d_{\mathbb{D}}}\left( {0,{F_\alpha }} \right)}}d\alpha } \le \frac{\pi }{2}\int_0^{ + \infty } {{\alpha ^{p - 1}}{\omega _{\mathbb{D}}}\left( {0,{F_\alpha }} \right)d\alpha }  <  + \infty.\]	
Conversely, suppose that for some $0<p<+\infty$,

\begin{equation}\label{hy}
\int_0^{ + \infty } {{\alpha ^{p - 1}}{e^{ - {d_{\mathbb{D}}}\left( {0,{F_\alpha }} \right)}}d\alpha }  <  + \infty.
\end{equation}
Without loss of generality, we set $\psi \left( 0 \right) = 0$. Let $dA$ denote the  Lebesgue measure on $\mathbb{D}$. For the Green function for $D$, set ${g_D}\left( {0,w} \right) = 0$ for $w \notin D$. By a change of variable and the conformal invariance of the Green function,

\begin{eqnarray}\label{rel}
\int_\mathbb{D} {{{\left| {\psi \left( z \right)} \right|}^{p - 2}}{{\left| {\psi '\left( z \right)} \right|}^2}\log \frac{1}{{\left| z \right|}}dA\left( z \right)} &=&\int_\mathbb{D} {{{\left| {\psi \left( z \right)} \right|}^{p - 2}}{{\left| {\psi '\left( z \right)} \right|}^2}{g_\mathbb{D}}\left( {0,z} \right)dA\left( z \right)} \nonumber \\
&=&\int_D {{{\left| w  \right|}^{p - 2}}{g_\mathbb{D}}\left( {0,{\psi ^{ - 1}}\left( w \right)} \right)dA\left( w  \right)} \nonumber \\
&=& \int_D {{{\left| w  \right|}^{p - 2}}{g_D}\left( {0,w } \right)dA\left( w \right)}        \nonumber \\
&=& \int_0^{ + \infty } {\int_0^{2\pi } {{\alpha ^{p - 2}}{g_D}\left( {0,\alpha {e^{i\theta }}} \right)\alpha d\theta d\alpha } } \nonumber \\
&=&\int_0^{ + \infty } {{\alpha ^{p - 1}}\left( {\int_0^{2\pi } {{g_D}\left( {0,\alpha {e^{i\theta }}} \right)d\theta } } \right)d\alpha }.
\end{eqnarray}
Applying elementary calculus, it is easily proved that there exist a positive constant $C$ and a point $x_0>0$ such that
\begin{equation}\label{cal}
\log \frac{{{e^x} + 1}}{{{e^x} - 1}} \le C {e^{-x}}
\end{equation}
for every $x \ge {x_0}$. Note that for $D$ unbounded and simply connected, ${d_D}\left( {0,\psi \left( {{F_\alpha }} \right)} \right) \to  + \infty $ as $\alpha \to  + \infty$ which also follows from the hypothesis (\ref{hy}). Therefore, by (\ref{cal}) and (\ref{gree}), we deduce that there exists an $\alpha_0>0$ such that for every $\alpha \ge \alpha_0$,

\[{g_D}\left( {0,\alpha {e^{i\theta }}} \right) \le C{e^{ - {d_D}\left( {0,\alpha {e^{i\theta }}} \right)}} \le C{e^{ - {d_D}\left( {0,\psi \left( {{F_\alpha }} \right)} \right)}} = C{e^{ - {d_\mathbb{D}}\left( {0,{F_\alpha }} \right)}}.\]
Integrating with respect to $\theta$, we get

\begin{equation}\label{in}
\int_0^{2\pi } {{g_D}\left( {0,\alpha {e^{i\theta }}} \right)d\theta }  \le C\int_0^{2\pi } {{e^{ - {d_\mathbb{D}}\left( {0,{F_\alpha }} \right)}}d\theta }  = 2\pi C{e^{ - {d_\mathbb{D}}\left( {0,{F_\alpha }} \right)}}
\end{equation}
for every $\alpha \ge \alpha_0$. So, by (\ref{rel}) and (\ref{in}), we infer that

\begin{equation}\label{fi}
\int_\mathbb{D} {{{\left| {\psi \left( z \right)} \right|}^{p - 2}}{{\left| {\psi '\left( z \right)} \right|}^2}\log \frac{1}{{\left| z \right|}}dA\left( z \right)} \le 2\pi C\int_{{\alpha _0}}^{ + \infty } {{\alpha ^{p - 1}}{e^{ - {d_\mathbb{D}}\left( {0,{F_\alpha }} \right)}}d\alpha }+ C',
\end{equation}
where 

\[C':=\int_0^{ \alpha_0 } {{\alpha ^{p - 1}}\left( {\int_0^{2\pi } {{g_D}\left( {0,\alpha {e^{i\theta }}} \right)d\theta } } \right)d\alpha }.\]
Finally, (\ref{hy}) and (\ref{fi}) give 

\[\int_\mathbb{D} {{{\left| {\psi \left( z \right)} \right|}^{p - 2}}{{\left| {\psi '\left( z \right)} \right|}^2}\log \frac{1}{{\left| z \right|}}dA\left( z \right)}<+\infty\]
and thus by \cite{Ya} we conclude that $\psi  \in {H^p\left( \mathbb{D} \right)}$.
\qed


\begin{bibdiv}
\begin{biblist}

\bib{Ahl}{book}{
title={Conformal Invariants: Topics in Geometric Function Theory},
author={L.V. Ahlfors},
date={1973},
publisher={McGraw-Hill},
address={New York}
}
\bib{Ba}{article}{
title={Univalent functions, Hardy spaces and spaces of Dirichlet type},
author={A. Baernstein and D. Girela and J. \'{A}. Pel\'{a}ez,},
journal={Illinois J. Math.},
volume={48},
date={2004},
pages={837--859}
}
\bib{Bea}{article}{
title={The hyperbolic metric and geometric function theory},
author={A.F. Beardon and D. Minda,},
journal={Quasiconformal mappings and their applications},
date={2007},
pages={9--56}
}
\bib{Du}{book}{
title={Theory of $H^p$ Spaces},
author={P.L. Duren},
date={1970},
publisher={Academic Press},
address={New York-London}
}
\bib{Gar}{book}{
title={Harmonic Measure},
author={J.B. Garnett and D.E. Marshall},
date={2005},
publisher={Cambridge University Press},
address={Cambridge}
}
\bib{Ha}{article}{
title={Hardy classes and ranges of functions},
author={L.J. Hansen},
journal={Michigan Math J.},
volume={17},
date={1970},
pages={235--248}
}
\bib{Ka}{article}{
title={On a relation between harmonic measure and hyperbolic distance on planar domains},
author={C. Karafyllia},
journal={Indiana Univ. Math. J. (to appear)}
}
\bib{Ra}{article}{
title={Univalent functions in Hardy, Bergman, Bloch and related spaces},
author={F. P\'{e}rez-Gonz\'{a}lez and J. R\"{a}tty\"{a}},
journal={J. d' Anal. Math.},
volume={105},
date={2008},
pages={125--148}
}
\bib{Co}{article}{
title={Geometric models, iteration and composition operators},
author={P. Poggi-Corradini},
journal={Ph.D. Thesis, University of Washington},
date={1996}
}
\bib{Ya}{article}{
title={Criteria for functions to be of Hardy class $H^p$},
author={S. Yamashita},
journal={Proc. Amer. Math. Soc.},
volume={75},
date={1979},
pages={69--72}
}

\end{biblist}
\end{bibdiv}
\end{document}